\documentclass{article}

\usepackage{amssymb}
\usepackage{graphicx}

\usepackage{amsmath}
\usepackage{amsfonts}
\usepackage{mathrsfs}
\usepackage[all]{xy}
\usepackage{latexsym}

\newtheorem{theorem}{Theorem}[section]

\newtheorem{lemma}[theorem]{Lemma}

\newtheorem{example}[theorem]{Example}

\newenvironment{pf}[1][\emph{Proof:}]{
\par\noindent{ #1}. }{\hfill\framebox(6,6)\par\medskip}

\begin{document}

\title{Cohomology of diagrams of algebras}
\author{Michael Robinson}
\date{\today}

\maketitle

\begin{abstract}
We consider cohomology of diagrams of algebras by Beck's approach,
using comonads. We then apply this theory to computing the
cohomology of $\Psi$-rings. Our main result is that there is a
spectral sequence connecting the cohomology of the diagram of an
algebra to the cohomology of the underlying algebra.
\end{abstract}

\section{Introduction}
Algebraic objects such as groups, Lie algebras, associative
algebras and commutative algebras have
cohomology theories defined in its own way. For example, group cohomology was
defined by Eilenberg-MacLane \cite{bib:group}, Lie algebra cohomology
was defined by Chevalley-Eilenberg \cite{bib:Liealg}, associative
algebra cohomology was defined by Hochschild \cite{bib:assoc alg},
and commutative algebra cohomology was defined by Andr\'{e}
\cite{bib:assoc alg1}, Quillen \cite{bib:assoc alg2}, and Barr
\cite{bib:assoc alg3}.  In the 1960's it was found that all of these
can be defined in one scheme. Here, we are going to use Barr-Beck's
approach, which is based on comonads. Let $T:\mathfrak{Sets}
\rightarrow\mathfrak{Sets}$ be a monad. Then one can consider the
category of $T$-algebras, $\mathfrak{Alg}(T)$, over the monad $T$. Let $A$ be a $T$-algebra
and $M$ be an $A$-module, which by definition is an abelian group object in $\mathfrak{Alg}(T)/A$. One
defines the cohomology of $A$ with coefficients in $M$, $H^{*}_{\mathbb{G}}(A,M)$, as follows.
There exists a pair of adjoint functors:
$$\xymatrix{\mathfrak{Sets}  \ar@<1ex>[rr]^{F} && \mathfrak{Alg}(T) \ar@<1ex>[ll]^{U}}$$
where $U$ is the forgetful functor. This adjoint pair of functors
yields a comonad $\mathbb{G}= FU: \mathfrak{Alg}(T) \rightarrow
\mathfrak{Alg}(T)$. One can take the comonad resolution
$\mathbb{G}_{*}(A)$. One can then apply the functor $Der(-,M)$. One
can then define a cochain complex by taking the alternating sum of
the induced maps, and then one can take the cohomology. This
situation is very general, one would like to apply this in the case
of $\lambda$-rings and $\Psi$-rings.

A $\psi$-rings is a commutative ring $R$ with identity $1$, together
with a series of homomorphisms $\psi^{n}:R \rightarrow R$, $n \geq
1$ such that $\forall$ $x \in R$ and integers $n,m \geq 1$ one has
$\psi^{1}(x) = x$ and $\psi^{n}(\psi^{m}(x)) = \psi^{nm}(x)$. So a
$\psi$-ring can be thought of as a diagram of an algebra with the monoid of the
natural numbers acting on $R$.

Recently, several authors have defined cohomology of diagrams of algebras.
For instance, the cohomology of diagrams of groups has been considered by
Cegarra \cite{bib:Cegarra}. Cohomology of diagrams of $\Pi$-algebras
has been considered by Blanc, Johnson, and Turner \cite{bib:blanc}.
Cohomology of diagrams of associative algebras was considered by
Gerstenhaber and Schack \cite{bib:Gersten1}. Cohomology of diagrams
of Lie algebras was considered by Gerstenhaber and Schack
\cite{bib:Gersten2}. All of these cohomologies are defined in its own way and it
is not clear how to make them as a particular case of one general theory.

The first aim is to use the Bar-Beck approach to unify the cohomology
of diagrams of algebras. Secondly, we would like to relate the
diagram cohomology to the cohomology of the algebra using a local to
global spectral sequence. Thirdly, we would like to apply the theory
to the category of commutative rings to study the cohomology of $\psi$-rings.

Our approach to defining the cohomology of diagrams of algebras can
be described as follows. First, fix a small category $I$. A diagram
of algebras is a functor $I \rightarrow \mathfrak{Alg}(T)$, where
$T$ is above a monad on sets. For appropriate $T$, one gets diagram
of groups, diagram of Lie algebras, commutative rings, etc.

One considers also the category $I_{0}$, which has the same objects
as $I$, but only the identity morphisms. The obvious inclusion
$I_{0} \subset I$ yields the functor $\mathfrak{Sets}^{I}
\rightarrow \mathfrak{Sets}^{I_{0}}$ which has left adjoint (the
left Kan extension). We also have the pair of adjoint functors
$\xymatrix{\mathfrak{Alg}(T)^{I} \ar@<0.5ex>[r] & \ar@<0.5ex>[l]
\mathfrak{Sets}^{I} }$ which comes from the adjoint pair
$\xymatrix{\mathfrak{Alg}(T) \ar@<0.5ex>[r] & \ar@<0.5ex>[l]
\mathfrak{Sets} }$. By gluing these diagrams, one gets another
adjoint pair $$\xymatrix{\mathfrak{Alg}(T)^{I} \ar@<0.5ex>[r] &
\ar@<0.5ex>[l] \mathfrak{Sets}^{I_{0}} }$$ We will prove that
$\mathfrak{Alg}(T)^{I}$ is monadic in $\mathfrak{Sets}^{I_{0}}$ and
the right cohomology theory of diagrams of algebras is one which is
associated to the corresponding comonad. These cohomology theories
are denoted by $H^{*}_{I}(A,M)$. We will prove that the comonad
cohomology is isomorphic to the ones considered in
\cite{bib:Cegarra} and \cite{bib:blanc} by choosing appropriate $T$.

The main technical element for studying $H^{*}_{I}(A,M)$ is the
local to global spectral sequence which can be described as follows.
Let $\textbf{A}:I \rightarrow \mathfrak{Alg}(T)$ be a diagram of
$T$-algebras and $M$ is an $\textbf{A}$-module. In particular, for
each $i \in I$ one has $M(i)$ an $A(i)$-module and for any arrow
$\alpha: i \rightarrow j$ one can consider $M(j)$ as an
$A(i)$-module where $A(i)$ acts on $M(j)$ via the algebra
homomorphism $A(\alpha):A(i) \rightarrow A(j)$. This allows us to
consider the cohomology $\mathcal{H}^{n}_{\mathbb{G}}(A,M)(\alpha)$
defined to be:
$$H^{n}_{\mathbb{G}}(A(i),\alpha^{*}M(j))$$
In this way we get a natural system \cite{bib:BW} on $I$.

Our main result claims that there exists a spectral sequence:
$$E^{pq}_{2} = H^{p}_{BW}(I,\mathcal{H}^{q}_{\mathbb{G}}(A,M)) \Rightarrow H^{*}_{I}(A,M).$$
where on the left hand side one uses the Baues-Wirsching cohomology
of a small category with coefficients in a natural system. This
spectral sequence is new in almost all of cases and gives
computational tools for diagram cohomology; even for diagrams of
groups, diagrams of associative algebras, diagrams of
$\Pi$-algebras, etc.

\section{Prerequisits}

\subsection{Baues Wirsching Cohomology}
For a small category $\mathcal{I}$. The \emph{category of
factorizations} in $\mathcal{I}$, denoted by $\mathcal{FI}$, is the
category with objects the morphisms $f,g,...$ in $\mathcal{I}$, and
morphisms $f\rightarrow g$ are pairs $(\alpha, \beta)$ of morphisms
in $\mathfrak{C}$ such that the following diagram commutes:

$$
\xymatrix{ B \ar[r]^{\alpha}  & B'
\\
A \ar[u]^{f} & A' \ar[l]^{\beta} \ar[u]_{g} }
$$
Composition in $\mathcal{FI}$ is given by
$(\alpha',\beta')(\alpha,\beta) = (\alpha' \alpha, \beta \beta')$. A
\emph{natural system} of abelian groups on $\mathcal{I}$ is a
functor from the category of factorizations to the category of
abelian groups:
$$D: \mathcal{FI} \rightarrow \mathfrak{Ab}$$
There are natural functors: $$\xymatrix{\mathcal{FI} \ar[r] &
\mathcal{I}^{op} \times \mathcal{I} \ar[r]  \ar[d] & \mathcal{I}
\ar[d] \\ & \mathcal{I}^{op} \ar[r] & \textbf{1}}$$ which allows one
to consider any functor or bifunctor on $\mathcal{I}$ as a natural
system on $\mathcal{I}$. Following Baues-Wirsching \cite{bib:BW},
one defines the cohomology $H_{BW}^{*}(\mathcal{I},D)$ of $\mathcal{I}$
with coefficients in the natural system $D$ as the cohomology of the
cochain complex $C^{*}_{BW}(\mathcal{I},D)$ given by:

$$
C_{BW}^{n}(\mathcal{I},D) = \prod_{ \alpha_{1}...\alpha_{n}:i_{n}
\rightarrow ... \rightarrow i_{0}} D(\alpha_{1}...\alpha_{n})
$$
and the coboundary map
$$d:C_{BW}^{n}(\mathcal{I},D) \rightarrow C_{BW}^{n+1}(\mathcal{I},D)$$
is given by:
\begin{align*}
(df)(\alpha_{1}...\alpha_{n+1}) &=
(\alpha_{1})_{*}f(\alpha_{2},...,\alpha_{n+1}) \\  &+
\sum_{j=1}^{n}(-1)^{i}f(\alpha_{1},...,\alpha_{j}\alpha_{j+1},...,\alpha_{n+1})
\\  &+ (-1)^{n+1}(\alpha_{n+1})^{*}f(\alpha_{1},...,\alpha_{n})
\end{align*}
We will need the following well-known lemma later:

\begin{lemma}\label{natsyslma}
Assume $i_{0} \in \mathcal{I}$ is an initial object, and
\\$F: \mathcal{I} \rightarrow Ab$ a functor. Then:
$$H_{BW}^{n}(\mathcal{I},F) = \left\{ \begin{array}{cc}
                     F(i_{0}) & n=0 \\
                       0 & n>0 \\
                      \end{array} \right.$$
\end{lemma}

\subsection{Base change}
Let $\mathfrak{C}$ be a category, and $X$ an object in $\mathfrak{C}$. An \emph{$X$-module} in $\mathfrak{C}$ is an abelian group
object in the category $\mathfrak{C}/X$.
$$X -mod := \mathrm{Ab}(\mathfrak{C}/X)$$

\begin{theorem}
Let $f: X \rightarrow Y$ be a morphism in $\mathfrak{C}$, then there
exists a base-change functor $f^{*}: Y-mod \rightarrow X-mod$ via
pullbacks.
\end{theorem}

\begin{pf}
The functor we are going to use is $f^{*}:\mathfrak{C}/Y \rightarrow
\mathfrak{C}/X $ given by pullbacks:

$$\xymatrix{f^{*}(M) \ar[r] \ar[d] & M \ar[d]^{p} \\ X \ar[r]^{f} & Y}$$
If $M \in Y-mod$ then $f^{*}(M)$ has a canonical $X$-module
structure. In set-theoretic notation:
$$f^{*}(M) = \{(x,m) | x \in X,\textrm{ } m \in M, \textrm{ } f(x) = p(m) \}$$
$$f^{*}(M) \times_{X} f^{*}(M) = \{(x,m, m') | x \in X,\textrm{ } m,m' \in M, \textrm{ } f(x) = p(m) = p(m')\}$$
$$f^{*}(M) \times_{X} f^{*}(M) \simeq f^{*}(M \times_{Y} M)$$
Consider the following commuting diagram.
$$\xymatrix{f^{*}(M \times_{Y} M) \ar[r] \ar[d] \ar@{..>}[ddr]^(0.3){\exists !} & M \times_{Y} M \ar[d] \ar[ddr]^{mult}\\ X \ar[r]_(0.35){f} \ar@{=}[ddr] & Y \ar@{=}[ddr] \\ & f^{*}(M) \ar[r] \ar[d] & M \ar[d]\\ & X \ar[r]_{f} & Y}$$
The unique morphism $f^{*}(mult):f^{*}(M \times_{Y} M) \rightarrow
f^{*}(M)$ exists by the universal property of pullbacks. The
isomorphism $f^{*}(M) \times_{X} f^{*}(M) \simeq f^{*}(M \times_{Y}
M)$ and this unique morphism yield multiplication:
$$f^{*}(mult):f^{*}(M) \times_{X} f^{*}(M) \rightarrow
f^{*}(M)$$ which gives an abelian group object structure on
$f^{*}(M)$.
\end{pf}

\subsection{Derivations}
For $M \in X$-mod, one defines a \emph{derivation} from $X$ to $M$
to be a morphism $d:X \rightarrow M$ which is a section of the
canonical morphism $M \rightarrow X$. Let $Der(X,M)$ denote the set
of derivations $d: X \rightarrow M$. We will require the following
useful theorem later.
\begin{theorem}\label{derttt}
If $X = \coprod_{\alpha \in I} X_{\alpha}$ and $M \in X$-mod, then
$Der(X,M) \cong \prod_{\alpha \in I} Der(X_{\alpha},M_{\alpha})$,
where $M_{\alpha}$ is an $X_{\alpha}$-module by the base-change
functor from the morphism $i_{\alpha}: X_{\alpha} \rightarrow X$.
\end{theorem}

\begin{pf}
From the definition of the coproduct one has a morphism $i_{\alpha}:
X_{\alpha} \rightarrow X$. Using this one gets $M_{\alpha} \in
X_{\alpha}$-mod via the following pullback diagram.
$$
\xymatrix{
 M  \ar[d]_{p} & M_{\alpha}  \ar[l]_{j_{\alpha}} \ar[d]^{p_{\alpha}}  \\
 X   & X_{\alpha} \ar[l]^{i_{\alpha}} }
$$
Let $f$ be a section of $p$, this means that $pf=id_{X}$. Consider
the following diagram:
$$
\xymatrix{ & & X_{\alpha} \ar@{..>}[ld]^{f_{\alpha}}
\ar[lld]_{fi_{\alpha}} \ar[ddl]^{id_{X_{\alpha}}} \\
 M  \ar@<1ex>[d]^{p}  & M_{\alpha}  \ar[l]_{j_{\alpha}} \ar[d]_{p_{\alpha}}  &  \\
 X  \ar@<1ex>[u]^{f}  &  X_{\alpha} \ar[l]^{i_{\alpha}} &  }
$$
The diagram commutes since $pfi_{\alpha} = id_{X}i_{\alpha} =
i_{\alpha}id_{X_{\alpha}}$. By the universal property of pullbacks $p_{\alpha}f_{\alpha} = id_{X_{\alpha}}$. So if $f$ is a section of
$p$ then $f_{\alpha}$ is a section of $p_{\alpha}$.

Conversely, let $f_{\alpha}$ be a section of $p_{\alpha}$, this
means that $p_{\alpha}f_{\alpha}=id_{X_{\alpha}}$. By the definition
of the coproduct there exists a unique $f$ such that the following
diagram commutes:
$$
\xymatrix{
 M  \\
 X  \ar@{..>}[u]^{f} &  X_{\alpha} \ar[l]^{i_{\alpha}} \ar[ul]_{j_{\alpha}f_{\alpha}} }
$$
This means that $fi_{\alpha}= j_{\alpha}f_{\alpha}$. Composing with $p$ on the left gives us that
$ pfi_{\alpha}= pj_{\alpha}f_{\alpha} = i_{\alpha}p_{\alpha}f_{\alpha} = i_{\alpha}id_{X_{\alpha}} = i_{\alpha}$
Thus the following diagram commutes:
$$
\xymatrix{ X \ar@<1ex>[r]^{id_{X}} \ar@<-1ex>[r]_{pf}  & X \\
X_{\alpha} \ar[u]^{i_{\alpha}} \ar[ur]_{i_{\alpha}}& & }$$
The universal property of the coproduct says that $pf =
id_{X}$. Hence $f$ is a section of $p$.
\end{pf}

\subsection{Coequalizers}
Let $f,g:a \rightarrow b$ be a pair in $\mathfrak{C}$, a
\emph{coequalizer} of $<f,g>$ is an arrow $u:b \rightarrow c$ such
that:
\begin{enumerate}
                                                                    \item $uf = ug$
                                                                    \item
if $h:b \rightarrow d$ with $hf = hg$, then there exists a unique
$h':c \rightarrow d$ such that $h = h'u$: $$\xymatrix{a
\ar@<0.5ex>[r]^{f} \ar@<-0.5ex>[r]_{g} & b \ar[dr]_{h} \ar[r]^{u} &
c \ar@{..>}[d]^{h'}
\\ && d}$$
                                                                  \end{enumerate}
$u$ is called a \emph{split coequalizer} of $f$ and $g$ if $u$ is a
coequalizer of $f$ and $g$ , and in addition there exists $s:c
\rightarrow b$ and $t:b \rightarrow a$ such that $us=1$ $ft=1$, and
$gt=su$:
$$\xymatrix{a
\ar@<0.5ex>[r]^{f} \ar@<-0.5ex>[r]_{g} & b \ar[r]^{u}
\ar@/^1pc/[l]^{t} & c \ar@/^1pc/[l]^{s} }$$

\subsection{Comonad Cohomology}
Let $\mathfrak{C}$ be a category, and $\mathbb{G}=$($G: \mathfrak{C}
\rightarrow \mathfrak{C}$, $\varepsilon: G \rightarrow Id$, $\delta:
G \rightarrow G^{2}$) be a comonad on $\mathfrak{C}$. For an object
$X$ in $\mathfrak{C}$, the comonad gives rise to a functorial
augmented simplicial object over $X$, which we denote by
$\mathbb{G}_{*}(X) \rightarrow X$. The object of $\mathbb{G}_{*}(X)$
in degree $n$ is $G^{n+1}(X)$, and the maps are $\varphi_{i} =
G^{i}\varepsilon G^{n-i}X: G^{n+1}(X) \rightarrow G^{n}(X)$, and
$\sigma_{i} = G^{i}\delta G^{n-i}X: G^{n+1}(X) \rightarrow
G^{n+2}(X)$ for $0 \leq i \leq n$. $X$ itself is in dimension $-1$
and the augmenting map is just $\varepsilon$.
$$
\xymatrix{ \ldots \ar@<2ex>[r] \ar@<-2ex>[r]^{\vdots} & G^{n}X
\ar@<2ex>[r] \ar@<-2ex>[r]^{\vdots} & G^{n-1}X \ar@<2ex>[r]
\ar@<-2ex>[r]^{\vdots} & \ldots \ar@<1ex>[r]^{G \varepsilon}
\ar@<-1ex>[r]_{\varepsilon G} & GX \ar[r]^{\varepsilon} & X}
$$
For any $M \in X-mod := \mathrm{Ab}$($\mathfrak{C}/X$), one can apply the functor
$Der(-,M)$ and take the alternating sum of the induced homomorphisms
to get the cochain complex whose cohomology is defined to be
$H_{\mathbb{G}}^{*}(X,M)$.

A morphism $f: X \rightarrow Y$ in $\mathfrak{C}$ is called a
\emph{$\mathbb{G}$-epimorphism} if the map
\\$Hom_{\mathfrak{C}}(G(Z),X) \rightarrow
Hom_{\mathfrak{C}}(G(Z),Y)$ is surjective for all $Z$. We require
the following useful lemma:

\begin{lemma}
$ \xymatrix{ GX \ar[r]^{\varepsilon_{X}} & X}$ is a
$\mathbb{G}$-epimorphism.
\end{lemma}

\begin{pf} For any map $h: GZ \rightarrow X$, we wish to find a map $f: GZ \rightarrow GX$ such that $f\varepsilon_{X} = h$.
We define $f$ via the following commutative diagram:

$$ \xymatrix{ G(GZ) \ar[r]^{G(h)} & G(X) \ar[r]^{\varepsilon_{X}} & X
\\
GZ \ar[u]^{\delta(Z)} \ar[ur]^{f}}$$

Now we need to check that $f\varepsilon_{X} = h$. By the naturality
of $\varepsilon$, the following diagram commutes:
$$\xymatrix{ GX \ar[r]^{\varepsilon_{X}} & X
\\
G(GZ) \ar[u]_{G(h)} \ar[r]_{\varepsilon_{GZ}} & GZ \ar[u]_{h}
\\
GZ \ar[u]_{\delta(Z)} \ar[ur]_{id_{GZ}} \ar@<2ex>@/^1pc/[uu]^{f}}$$
So $\varepsilon_{X}$ is a $\mathbb{G}$-epimorphism.
\end{pf}

An object $P$ of $\mathfrak{C}$ is called
\emph{$\mathbb{G}$-projective} if for any $\mathbb{G}$-epimorphism
$f:X \rightarrow Y$, the map $Hom_{\mathfrak{C}}(P,X) \rightarrow
Hom_{\mathfrak{C}}(P,Y)$ is surjective. Later we will require the
following lemmas:

\begin{example}
For all $Z$, $G(Z)$ is $\mathbb{G}$-projective.
\end{example}

\begin{lemma}
The coproduct of $\mathbb{G}$-projective objects is
$\mathbb{G}$-projective.
\end{lemma}

\begin{pf}
Let $P = \coprod_{i} P_{i}$ where $P_{i}$ is $\mathbb{G}$-projective
for all $i$. For a map \\$f: X \rightarrow Y$, one applies the
functors $Hom_{\mathfrak{C}}(P,-)$ and $Hom_{\mathfrak{C}}(P_{i},-)$
to get the maps $f_{*}: Hom_{\mathfrak{C}}(P,X) \rightarrow
Hom_{\mathfrak{C}}(P,Y)$, $f_{i*}: Hom_{\mathfrak{C}}(P_{i},X)
\rightarrow Hom_{\mathfrak{C}}(P_{i},Y)$. If $f$ is a
$\mathbb{G}$-epimorphism then $f_{i*}$ is surjective for all $i$.
Using the well-known lemma $Hom_{\mathfrak{C}}(\coprod_{i}P_{i},Z)
\cong \prod_{i}Hom_{\mathfrak{C}}(P_{i},Z)$ one gets that if $f$ is
a $\mathbb{G}$-epimorphism then $f_{*} \cong \prod_{i} f_{i*}$ is
surjective. Hence $P$ is $\mathbb{G}$-projective if $P_{i}$ is
$\mathbb{G}$-projective for all $i$.\end{pf}

\begin{lemma}
An object $P$ is $\mathbb{G}$-projective if and only if $P$ is a
retract of an object of the form $G(Z)$.
\end{lemma}

\begin{pf}
A retract of a surjective map is surjective, so it is sufficient to
consider the case $P=G(Z)$, which is obvious from the definition of
$\mathbb{G}$-epimorphism.
\end{pf}

\begin{lemma}\label{exacttt}
$H^{p}_{\mathbb{G}}(X,M)=0$, for $p>0$ provided $X$ is
$\mathbb{G}$-projective.
\end{lemma}

\begin{pf}
From the previous lemma, it is sufficient to check this in the case
where $X = G(Z)$. In this case it is possible to construct a
contracting homotopy. There are maps $s_{n}: G^{n+2} \rightarrow
G^{n+3}$ for $n \geq -1$ such that $\epsilon s_{-1} = id$,
$\varphi_{n+1}s_{n} = id$, $\varphi_{0}s_{0} = s_{-1}\epsilon$, and
$\varphi_{i}s_{n} = s_{n-1}\varphi_{i}$ for all $0 \leq i \leq n$.
$$s_{n} = G^{n+1}\delta$$
It follows that $H^{p}_{\mathbb{G}}(G(Z),M)=0$, for $p>0$.
\end{pf}

\begin{lemma}
$H^{0}_{\mathbb{G}}(X,M)= Der(X,M)$ for all $X$.
\end{lemma}

\section{$T$-algebras} Start with an adjunction $\mathfrak{Sets}
\rightharpoonup \mathfrak{C} $, and construct the monad
$T:\mathfrak{Sets} \rightarrow\mathfrak{Sets}$ defined in
$\mathfrak{C}$. Then one can consider the category of $T$-algebras,
$\mathfrak{Alg}(T)$, over the monad $T$. There exists an adjoint
pair of functors:
$$\xymatrix{\mathfrak{Sets}  \ar@<1ex>[rr]^{F} && \mathfrak{Alg}(T) \ar@<1ex>[ll]^{U}}$$
where $U$ is the forgetful functor. This adjoint pair of functors
yields a comonad $\mathbb{G}= FU: \mathfrak{Alg}(T) \rightarrow
\mathfrak{Alg}(T)$. There exists a unique functor $K: \mathfrak{C}
\rightarrow \mathfrak{Alg}(T)$. If $K$ is an equivalence of
categories, then one says that $\mathbb{G}$ is \textit{monadic}.

\begin{theorem}[Beck's Theorem \cite{bib:Mac}]\label{Beck theorem}
The following are equivalent:
\begin{enumerate}
  \item The comparison functor $K: \mathfrak{C}
\rightarrow \mathfrak{Alg}(T)$ is an equivalence of categories.
  \item If $f,g$ is any parallel pair in $\mathfrak{C}$ for which $Uf$, $Ug$
has a split coequalizer, then $\mathfrak{C}$ has a coequalizer for
$f,g,$ and $U$ preserves and reflects coequalizers for these pairs.
\end{enumerate}
\end{theorem}
Fix a small category $I$, one also considers the category $I_{0}$,
which has the same objects as $I$, but only the identity morphisms.
The obvious inclusion $I_{0} \subset I$ yields the functor
$\mathfrak{Sets}^{I} \rightarrow \mathfrak{Sets}^{I_{0}}$ which has
left adjoint (the left Kan extension). Let $\mathcal{F}: I_{0}
\rightarrow \mathfrak{Sets}$ be a functor, then the left Kan
extension $Lan(\mathcal{F}): I \rightarrow \mathfrak{Sets}$ is given
by $Lan(\mathcal{F})(i) = \coprod_{x \rightarrow i}F(x)$. We also
have the pair of adjoint functors $\xymatrix{\mathfrak{Alg}(T)^{I}
\ar@<0.5ex>[r] & \ar@<0.5ex>[l] \mathfrak{Sets}^{I} }$ which comes
from the adjoint pair $\xymatrix{\mathfrak{Alg}(T) \ar@<0.5ex>[r] &
\ar@<0.5ex>[l] \mathfrak{Sets} }$. By gluing these diagrams, one
gets another adjoint pair $$\xymatrix{\mathfrak{Alg}(T)^{I}
\ar@<0.5ex>[r]^{U_{I}} & \ar@<0.5ex>[l]^{F_{I}}
\mathfrak{Sets}^{I_{0}} }$$This adjoint pair of functors yields a
comonad $\mathbb{G}_{I}= F_{I}U_{I}: \mathfrak{Alg}(T)^{I}
\rightarrow \mathfrak{Alg}(T)^{I}$. If $A: I \rightarrow
\mathfrak{Alg}(T)$, then $$\mathbb{G}_{I}(A)(i) = \coprod_{x
\rightarrow i}\mathbb{G}(x)$$

\begin{lemma}
If $\mathbb{G}$ is monadic, then $\mathbb{G}_{I}$ is monadic.
\end{lemma}

\begin{pf}
Assume $\mathbb{G}$ is monadic and consider a parallel pair $f,g$ in
$\mathfrak{C}^{I}$:$$\xymatrix{F \ar@<0.5ex>[r]^{f}
\ar@<-0.5ex>[r]_{g} & T}$$ If there is a split coequalizer in $X$ as
follows:
$$\xymatrix{UF \ar@<0.5ex>[r]^{Uf} \ar@<-0.5ex>[r]_{Ug} & UT
\ar[r]^{q} & W}$$ then by theorem \ref{Beck theorem} for each $i \in
I$, one has that $q(i)$ is a coequalizer of the following:
$$\xymatrix{UF(i) \ar@<0.5ex>[r]^{Uf(i)} \ar@<-0.5ex>[r]_{Ug(i)} &
UT(i)}$$ Hence by theorem \ref{Beck theorem} there exists $Z(i)$,
$h(i)$ in $\mathfrak{C}$ with :$$\xymatrix{F(i)
\ar@<0.5ex>[r]^{f(i)} \ar@<-0.5ex>[r]_{g(i)} & T(i) \ar[r]^{h(i)} &
Z(i)}$$ such that $UZ(i) = W(i)$ and $Uh(i) = q(i)$. In fact $Z$ is
a functor $I \rightarrow \mathfrak{C}$. For $\alpha:i \rightarrow j$
one considers the following commuting diagram for $(\alpha: i
\rightarrow j) \in I$:$$\xymatrix{F(i) \ar[d]_{F(\alpha)}
\ar@<0.5ex>[r]^{f(i)}
\ar@<-0.5ex>[r]_{g(i)} & T(i) \ar[d]^{T(\alpha)} \ar[r]^{h(i)} & Z(i)\ar@{..>}[d]^{\exists !} \\
F(j) \ar@<0.5ex>[r]^{f(j)} \ar@<-0.5ex>[r]_{g(j)} & T(j)
\ar[r]^{h(j)} & Z(j)}$$ Since the coequalizer is universal, this
means that there exists a unique map $Z(\alpha): Z(i) \rightarrow
Z(j)$ which commutes with the above diagram. One can check that $Z$
is indeed a functor.
\end{pf}

\begin{lemma}\label{derlma}
For all objects $Z$ in $\mathfrak{C}$, and for $A \in G(Z)-mod$, one
has
$$Der(G(Z),A) = \{ s \in \mathfrak{Sets}(UZ,UA) | U(\pi)s = j_{UZ}:U(Z) \rightarrow UFU(Z)\}$$ where $\pi$ is the canonical morphism $\pi: A \rightarrow G(Z)$
\end{lemma}

\begin{pf}
From the definition of $F,U$ being an adjoint pair, one gets that:
$$Hom_{\mathfrak{Alg(\mathbb{T})}}(F(X),Y) =
Hom_{\mathfrak{Sets}}(X,U(Y))$$ Setting $X= U(Z)$ and $Y = FU(Z)$ we
get the following:$$Hom_{\mathfrak{Alg(\mathbb{T})}}(FU(Z),FU(Z)) =
Hom_{\mathfrak{Sets}}(U(Z),UFU(Z))$$ From this it can be shown that
$$Der(G(Z),A) = \{ s \in \mathfrak{Sets}(UZ,UA) | U(\pi)s =
j_{UZ}:U(Z) \rightarrow UFU(Z)\}$$
\end{pf}

\section{Cohomology of diagrams of algebras}
In this section, let $\mathfrak{C}$ denote the category of sets, and
$I$ denote a small category. $\mathfrak{C}$ is a category with
limits. We require the following useful theorem.

\begin{theorem}\label{natsysttt}
Let $A: I \rightarrow \mathfrak{C}$ be a functor, and $M \in A$-mod
$:= Ab(\mathfrak{C}^{I}/A)$ where $\mathfrak{C}^{I}$ is the category
of functors. If ($\alpha : i \rightarrow j$) $\in I$, then $M(j) \in
A(j)-mod$ and
$$ \mathfrak{Der}(A,M)(\alpha) = Der(A(i), \alpha^{*}M(j)) $$
defines a natural system on $I$.
\end{theorem}

\begin{pf}
Start by fixing $A$ and $M$, then let $D(\alpha)$ denote $
\mathfrak{Der}(A,M)(\alpha)$. Let $\gamma, \alpha, \beta \in I$ such
that:
$$
\xymatrix{ i' \ar[r]^{\gamma} & i \ar[r]^{\alpha} & j \ar[r]^{\beta}
& j' }
$$
We are going to show that:
$$
\xymatrix{ D(\alpha \gamma) & D(\alpha) \ar[l]_-{\gamma^{*}}
\ar[r]^-{\beta^{*}} & D(\beta \alpha) }
$$
Let $s \in D(\alpha)$, then the following diagram commutes with $ps
= id_{A(i)}$, and $\alpha^{*}M(j)$ is a pullback; $\alpha^{*}M(j)
\in A(i)-mod.$

$$
\xymatrix{\alpha^{*}M(j) \ar[rr]  \ar@<1ex>[dd]^{p} && M(j) \ar[dd] \\
 & &\\
A(i) \ar@<1ex>[uu]^{s} \ar[rr]^{A(\alpha)} && A(j) }$$

Consider the following commuting diagram:

$$
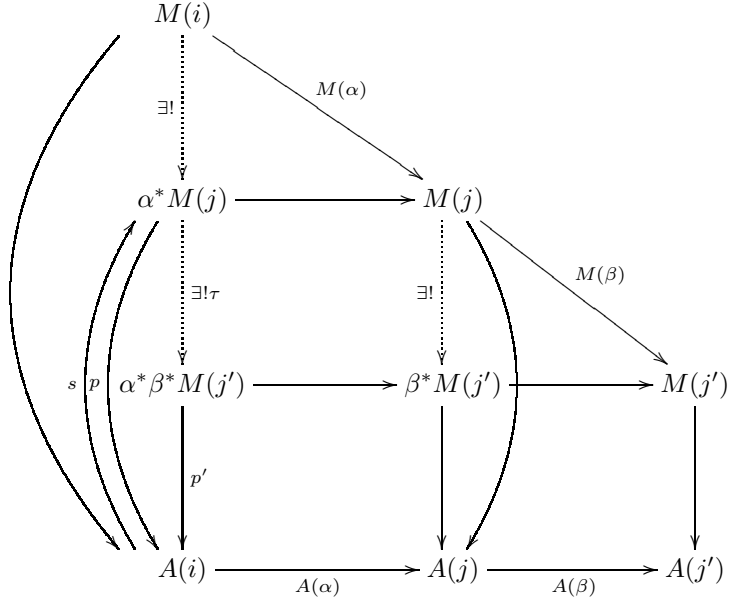
\xymatrix{M(i) \ar@<-4ex>@/_4pc/[dddddd] \ar@{..>}[dd]_{\exists !} \ar[ddrr]^{M(\alpha)} \\
\\
 \alpha^{*}M(j) \ar@<-1ex>@/_2pc/[dddd]_{p} \ar@{..>}[dd]^{\exists ! \tau} \ar[rr] && M(j) \ar@/^2pc/[dddd] \ar@<-1ex>@{..>}[dd]_{\exists !} \ar[ddrr]^{M(\beta)} \\
&  \\
\alpha^{*}\beta^{*}M(j') \ar[rr] \ar[dd]^{p'}  && \beta^{*}M(j') \ar@<-1ex>[dd] \ar[rr] && M(j') \ar[dd]\\
&  & & && \\
A(i) \ar@<3ex>@/^2pc/[uuuu]^{s} \ar[rr]_{A(\alpha)} && A(j)
\ar[rr]_{A(\beta)} && A(j')}
$$
Let $s': A(i) \rightarrow \alpha^{*}\beta^{*}M(j')$ be the map $s' = \tau s$.
If we let $s \in D(\alpha)$, this means that $ps = id_{A(i)}$. Hence $$p' \tau s = p s
= id_{A(i)}$$ Hence $s' \in Der(A(i), \alpha^{*}\beta^{*}M(j')) =
Der(A(i), (\beta\alpha)^{*}M(j'))$.

Consider the following commutative diagram, with $s$ a section of
$p$.
$$
\xymatrix{
(\alpha \gamma)^{*}M(j) \ar[rr]  \ar@<1ex>[dd]^{p'} && \alpha^{*}M(j) \ar[rr]  \ar@<1ex>[dd]^{p} && M(j) \ar[dd] \\
&  & &  &\\
A(i')  \ar[rr]^{A(\gamma)} && A(i) \ar@<1ex>[uu]^{s}
\ar[rr]^{A(\alpha)} && A(j) }$$ There exists a unique $s': A(i')
\rightarrow (\alpha \gamma)^{*}M(j)$ which is a section of $p'$
which would make the above diagram still commute. Therefor $s' \in
Der(A(i'), (\alpha\gamma)^{*}M(j))$.
\end{pf}

Let $\mathbb{G}$ be a comonad in $\mathfrak{C}$. Let $A: I
\rightarrow \mathfrak{C}$ be a functor, and $M \in A$-mod. Then we
can construct the following bicomplex denoted by $C^{*,*}(I,A,M)$:

$$
C^{p,q}(I,A,M) = \prod_{i_{0}\rightarrow ... \rightarrow
i_{p}}Der(G^{q+1}(A(i)),M(k))
$$
The map $C^{p,q}(I,A,M) \rightarrow C^{p+1,q}(I,A,M)$ is the map in
the Baues-Wirsching cochain complex, and the map $C^{p,q}(I,A,M)
\rightarrow C^{p,q+1}(I,A,M)$ is the coproduct of maps in the
comonad cochain complex.

$$
\xymatrix{ \vdots & \vdots & \vdots
\\\prod_{i}Der(G^{3}(A(i)),M(i)) \ar[u]\ar[r] &
\prod_{i\rightarrow j}Der(G^{3}(A(i)),M(j)) \ar[r]\ar[u] &
\prod_{i\rightarrow j\rightarrow k}Der(G^{3}(A(i)),M(k)) \ar[r]
\ar[u]& \ldots \\\prod_{i}Der(G^{2}(A(i)),M(i)) \ar[r] \ar[u]&
\prod_{i\rightarrow j}Der(G^{2}(A(i)),M(j)) \ar[r]\ar[u] &
\prod_{i\rightarrow j\rightarrow k}Der(G^{2}(A(i)),M(k)) \ar[r]
\ar[u]& \ldots
\\\prod_{i}Der(G(A(i)),M(i)) \ar[r] \ar[u] & \prod_{i\rightarrow
j}Der(G(A(i)),M(j)) \ar[r]  \ar[u] & \prod_{i\rightarrow
j\rightarrow k}Der(G(A(i)),M(k)) \ar[u] \ar[r] & \ldots}
$$
We let \emph{$H^{*}(I,A,M)$} denote the cohomology of the total
complex of  $C^{*,*}(I,A,M)$.

We will need the following useful lemmas:

\begin{lemma}\label{GIprojttt}
If $A$ is $\mathbb{G}_{I}$-projective, then $A(i)$ is $\mathbb{G}$-projective for all
$i \in I$.
\end{lemma}

\begin{pf}
Consider $A = G_{I}(Z): I \rightarrow \mathfrak{C}$ where
$G_{I}(A)(i) = \coprod _{x \rightarrow i} G(A(x))$. Since $G(A(x))$
is $\mathbb{G}$-projective, it follows that $\coprod _{x \rightarrow i}
G(A(x))$ is $\mathbb{G}$-projective for all $i \in I$.
\end{pf}

\begin{lemma}\label{lemmaHIAM}
$H^{0}(I,A,M) = \mathfrak{Der}(A,M)$, furthermore, if $A$ is
$\mathbb{G}_{I}$-projective then $H^{n}(I,A,M) = 0$ for $n > 0$.
\end{lemma}

\begin{pf}
It is sufficient to consider the case when $A= G_{I}(Z)$. When $A =
G_{I}(Z)$, it is known that $A$ is $\mathbb{G}_{I}$-projective. By
theorems \ref{GIprojttt} and \ref{exacttt}, one gets that the
vertical columns in our bicomplex are exact except in dimension 0.
There is a well known lemma for bicomplexes which tells us the
cohomology of the total complex is isomorphic to the cohomology of
the following chain complex:

$$\xymatrix{\prod_{i}Der(A(i),M(i)) \ar[r] & \prod_{i \rightarrow j}Der(A(i),M(j)) \ar[r] & \ldots}$$
It is known that the cohomology of this cochain complex is just
$H_{BW}^{*}(I,\mathfrak{Der}(A,M))$.

To prove the first statement it is enough to show that
$$0 \rightarrow \mathfrak{Der}(A,M) \rightarrow \prod_{i}Der(A(i),M(i)) \rightarrow \prod_{\alpha:i \rightarrow j}Der(A(i),M(j))$$
is exact. Let $\psi \in \prod_{i}Der(A(i),M(i))$ and $(\alpha:i \rightarrow j) \in I$, then $d \psi (\alpha : i \rightarrow j) = \alpha_{*}\psi(i) - \alpha^{*}\psi(j).$ Therefore $d\psi(\alpha: i \rightarrow j) =0$
if and only if $\alpha_{*} \psi(i) = \alpha^{*}\psi(j)$. However $\alpha_{*} \psi(i) = \alpha^{*}\psi(j)$ if and
only if $M(\alpha)\psi(i) = \psi(j)A(\alpha),$ i.e. the following
diagram commutes:
$$\xymatrix{A(i) \ar[d]^{A(\alpha)} \ar[r]^{\psi(i)} & M(i) \ar[d]^{M(\alpha)} \\ A(j) \ar[r]_{\psi(j)} & M(j)}$$
Hence $\psi \in \mathfrak{Der}(A,M)$. This tells us that the sequence above is exact. Hence $H^{0}(I,A,M) = \mathfrak{Der}(A,M)$.

To prove the second statement, let us consider \begin{align*}
D(\alpha: i\rightarrow j) :&= \mathfrak{Der}(A(i),\alpha^{*}M(j))
\\ & = Der(\coprod_{\beta:y \rightarrow i}GZ(y), \alpha^{*}M(j))\\
& = \prod_{\beta:y \rightarrow i}Der(GZ(y),
\beta^{*}\alpha^{*}M(j)), \textrm{ by lemma \ref{derttt}}
\end{align*}

Define $D_{y}$ for a fixed object $y \in I$ to be a natural system
on $I$ (using theorem \ref{natsysttt}) given by:
$$D_{y}(\alpha:i \rightarrow j) = \prod_{\beta:y \rightarrow i}Der(GZ(y), \beta^{*}\alpha^{*}M(j))$$

So one has that:
$$D(i\rightarrow j) = \prod_{y}D_{y}(i\rightarrow j)$$
Hence,
$$H^{*}_{BW}(I,D) = \prod_{y \in I}H^{*}_{BW}(I,D_{y})$$

Now consider the cochain complex $C^{*}_{BW}(I,D_{y})$:

\begin{align*}
C^{*}_{BW}(I,D_{y}) &= \xymatrix{\prod _{i}D_{y}(i\rightarrow i) \ar[r] &
\prod_{\alpha: i \rightarrow j}D_{y}(i\rightarrow j) \ar[r] &
\ldots} \\ &= \xymatrix{\prod_{i} \prod_{\beta:y \rightarrow i}Der(GZ(y), \beta^{*}M(i)) \ar[r] & \prod_{ \alpha : i \rightarrow j} \prod_{\beta:y \rightarrow i}Der(GZ(y),\beta^{*}\alpha^{*}M(j)) \ar[r] & \ldots}
\end{align*}
$UZ(y)$ forms a basis of the free object $GZ(y)$, applying lemma
\ref{derlma}, one can rewrite the cochain complex as:

$$C^{*}_{BW}(I,D_{y}) = \xymatrix{\prod_{y\rightarrow i} \prod_{m
\in UZ(y)}A_{\beta j(m)} \ar[r]& \prod_{\alpha:i \rightarrow
j}\prod_{\beta:y \rightarrow i}\prod_{m \in UZ(y)}A_{\alpha\beta
j(m)} \ar[r] & \ldots } $$ where $A_{\beta j(m)} = \textrm{preimage
of } \beta\gamma(m) \textrm{ in the projection } M(j) \rightarrow
GZ(j)$. This allows us to rewrite the cochain complex as
$$C^{*}_{BW}(I,D_{y}) = \prod_{m \in UZ(y)}C^{*}_{BW}(y/I, F_{m})$$
where $F_{m}: y/I \rightarrow Ab$ is a functor defined by
$F_{m}(\beta:y\rightarrow i) = A_{\beta j(m)}$

Since the category $y/\mathfrak{C}$ contains an initial object
($id_{y}:y \rightarrow y$), so by lemma \ref{natsyslma} the
cohomology vanishes in positive dimensions.
\end{pf}

\begin{theorem}
$H^{*}_{G_{I}}(A,M) = H^{*}(I,A,M)$
\end{theorem}

\begin{pf}
We are going to show that: $$H^{*}(C^{\bullet}(I,A,M)) \simeq
H^{*}(Tot(C^{\bullet}(I,G_{I}(A),M))) \simeq
H^{*}(C^{\bullet}_{G_{I}}(A,M))$$

Start by considering $C^{\bullet}(I,G_{I}^{p}(A),M)$. Since
$G_{I}^{p}(A)$ is $\mathbb{G}_{I}$-projective, it follows that
$$H^{n}(Tot(C^{\bullet}(I,G_{I}(A),M))) = \left\{
                                           \begin{array}{ll}
                                             \mathfrak{Der}(G_{I}^{p}(A),M), & \hbox{n = 0;} \\
                                             0, & \hbox{otherwise.}
                                           \end{array}
                                         \right.
$$
Hence by lemma \ref{lemmaHIAM} $H^{*}(C^{\bullet}(I,A,M)) \simeq
H^{*}(Total(C^{\bullet}(I,G_{I}(A),M)))$.

Now let us consider $C^{p,q}(I,G_{I}(A),M) = \prod_{x_{0}\rightarrow
... \rightarrow
x_{p}}\mathfrak{Der}(G^{q+1}(G^{*}_{I}(A)(x_{0})),M(x_{p}))$. One
has that $\mathbb{G}^{*}_{I}(A) \rightarrow A$, which is an
augmented simplicial object and one can apply the functor $U_{I}$ to
get: $U_{I}G^{*}_{I}(A) \rightarrow U_{I}(A)$ which is contractible
in $(\mathfrak{Sets}^{I_{0}})$. Then one can apply the functor
$F_{I}$ to get $G_{I} \mathbb{G}^{*}_{I}(A) \rightarrow G_{I}(A)$
which is contractible in $\mathfrak{Alg}(T)$. Hence $G^{q+1}_{I}
\mathbb{G}^{*}_{I}(A) \rightarrow G^{q+1}_{I}(A)$ is contractible in
$\mathfrak{Alg}(T)$. Applying the functor $\mathfrak{Der}(-,M)$, one
gets a contractible cosimplicial abelian group. Hence
$H^{*}(Tot(C^{\bullet}(I,G_{I}(A),M))) \simeq
H^{*}(C^{\bullet}_{G_{I}}(A,M))$.
\end{pf}

Now one has both a global cohomology, $H^{*}_{\mathbb{G}_{I}}(A,M)$,
and a local cohomology, $H^{*}(A(i),M(i))$. One can ask how these
two are related; the answer is given by by the local to global
spectral sequence:
$$E^{pq}_{2} = H^{p}_{BW}(I,\mathcal{H}^{q}(A,M)) \Rightarrow H^{p+q}_{\mathbb{G}_{I}}(A,M).$$
where $\mathcal{H}^{q}(A,M)$ is a natural system on $I$ whose value
on $(\alpha: i \rightarrow j)$ is given by $H^{q}(A(i),M(j))$.

\section{Applications}
\subsection{$\Psi$-rings} As an example of the general theory, one
can consider $\Psi$-rings. A \emph{$\Psi$-ring} is a commutative
ring $R$ with identity $1$, with a sequence of ring homomorphisms
$\Psi^{n}: R \rightarrow R$, $n \geq 1$ satisfying $\forall x \in
R$, and integers $n,m \geq 1$.
\begin{enumerate}
    \item $\Psi^{1}(x) = x$
    \item $\Psi^{n}( \Psi^{m}(x)) = \Psi^{nm}(x)$
\end{enumerate}

So to know $\Psi$-rings it is sufficient to know $\Psi^{p}$, for $p$
prime such that for all primes $p,q$, $\Psi^{p} \Psi^{q} = \Psi^{q}
\Psi^{p}$.

If $R$ is a $\Psi$-ring, $M$ is a \emph{$\Psi$-module} if $M$ is an
$R$-module together with a sequence of homomorphisms $\Psi^{n}:M
\rightarrow M$ such that for all $m \in M$, $r \in R$, $l,n \geq 1$:
\begin{enumerate}
    \item $\Psi^{1}(m) = m $
    \item $\Psi^{n}(rm) = \Psi^{n}(r)\Psi^{n}(m) = \Psi^{n}(m) \Psi^{n}(r)$
    \item $\Psi^{n}(\Psi^{l}(m))  = \Psi^{nl}(m)$
\end{enumerate}
We let $R$-mod$_{\Psi}$ denote the category of all $\Psi$-modules
over $R$.

Let $R$ be a $\Psi$-ring, $M \in R$-mod$_{\Psi}$ then the semidirect
product of the underlying ring and module, $R \rtimes M$, together
with maps: $\Psi^{i}: R \rtimes M \rightarrow R \rtimes M$ for $i
\geq 1$ given by: $$ \Psi^{i} (r,m) = (\Psi^{i}(r), \Psi^{i}(m))$$
is a $\Psi$-ring. We call this the \emph{semi-direct product} of $R$
and $M$, denoted by $R \rtimes_{\Psi} M$

We define a \emph{$\Psi$-derivation} is a $\Psi$-module homomorphism
$d:R \rightarrow M$ such that $\forall r, r' \in R, \forall n \geq
1$:
\begin{displaymath}
d(r r') = r d(r') + d(r)r'
\end{displaymath}
\begin{displaymath}
\Psi^{n}(d(r)) = d(\Psi^{n}(r))
\end{displaymath}
We let $Der_{\Psi}(R,M)$ denote the set of all $\Psi$-derivations
$d:R \rightarrow M$. One would expect the following theorem:

\begin{theorem}
There is a one-to-one correspondence between the sections of $
\xymatrix{ R \rtimes_{\Psi} M \ar[r]^-{\pi} & R }$ and the
$\Psi$-derivations $d:R \rightarrow M$
\end{theorem}

\begin{pf}[Proof of theorem]
Assume we have a section of $\pi$, then:$$ \xymatrix{ R
\rtimes_{\Psi} M \ar@<1ex>[r]^-{\pi} & R \ar@<1ex>[l]^-{\sigma} }$$
$\sigma \pi = id_{R}$, so $\sigma(x) = (x,d(x))$ for some $d:R
\rightarrow M$.
$$d(x+y) = d(x) + d(y), \qquad d(xy) = d(x)y + x d(y)$$ follow from
$\sigma$ being a homomorphism of $\Psi$-rings. $\sigma$ preserves
the $\Psi$-ring structure, meaning that $\Psi^{i}\sigma(x) =
\sigma\Psi^{i}(x)$.
$$\Psi^{i}\sigma(x) = \Psi^{i}(x, d(x)) = (\Psi^{i}(x),\Psi^{i}(d(x))) $$  $$\sigma\Psi^{i}(x) = (\Psi^{i}(x),d(\Psi^{i}(x)))
$$ Hence $\Psi^{i}\sigma(x) = \sigma\Psi^{i}(x)$ if and only if
$\Psi^{i}d(x) = d\Psi^{i}(x)$. This tells us that if $\sigma$ is a
section of $\pi$, then we have a $\Psi$-derivation $d$.

Conversely, if we have a $\Psi$-derivation $d:R \rightarrow M$, then
$\sigma(x) = (x,d(x))$ is a section of $\pi$.
\end{pf}

We now construct the free $\Psi$-ring on one generator $a$. Let $A$
be the free commutative ring generated by $a_{0}, a_{1}, a_{2},...$.
Since there are countably infinitely many primes, it is possible to
label them with the natural numbers. Set $a_{0} = a$, and $a_{i} =
\Psi^{p}(a)$, where $p$ is the $i^{th}$ prime, for $i \in
\mathbb{N}$. Then $A$ is a $\Psi$-ring.

More generally, we can construct a free $\Psi$-ring on generators
$a, b, c,... n$. We let $R$ be the free commutative ring generated
by $a_{0}, a_{1},...,b_{0},
b_{1},...,c_{0},c_{1},...,n_{0},n_{1},...$. Set $a_{0} = a$, $b_{0}
= b$, $c_{0} = c$, $\ldots$, $n_{0} = n$, and $a_{i} = \Psi^{p}(a)$,
$b_{i} = \Psi^{p}(b)$,$c_{i} = \Psi^{p}(c)$, $\ldots$, $n_{i} =
\Psi^{p}(n)$ where $p$ is the $i^{th}$ prime, for $i \geq 1$. Then
$R$ is a $\Psi$-ring.

There is a forgetful functor $U: \Psi\mathfrak{-rings} \rightarrow
\mathfrak{Sets}$ from the category of \\$\Psi$-rings to the category
of sets. This has the left adjoint $F: \mathfrak{Sets} \rightarrow
\Psi\mathfrak{-rings}$, where $F(S)$ is the free $\Psi$-ring
generated by $S \in \mathfrak{Sets}$. Hence there is an adjoint pair
of functors:
$$\xymatrix{\mathfrak{Sets}  \ar@<1ex>[rr]^{F} && \Psi\mathfrak{-rings} \ar@<1ex>[ll]^{U}}$$
where $U$ is the forgetful functor, and $F$ is the free functor. The
adjoint pair of functors yields a comonad $\mathbb{G} = FU:
\Psi\mathfrak{-rings} \rightarrow \Psi\mathfrak{-rings}$ which is
monadic.

Let $\mathbb{N}^{mult}$ denote the multiplicative monoid of the
natural numbers, and let $I$ denote the category with one object
associated to $\mathbb{N}^{mult}$. Then one can consider
$\Psi$-rings as diagrams of algebras  being functors
 from $I$ to the category of commutative rings, $\mathfrak{Com.rings}$.
So $\Psi$-rings are diagrams of algebras with $\mathbb{N}^{mult}$
acting on $R$ a commutative ring with identity. Hence we can use the
theory which we developed in the previous section.

It is well known that there is an adjoint pair of functors:
$$\xymatrix{\mathfrak{Sets}  \ar@<1ex>[rr]^{F} && \mathfrak{Com.rings} \ar@<1ex>[ll]^{U}}$$
This gives rise to a comonad $\mathbb{G} = FU: \mathfrak{Com.rings}
\rightarrow \mathfrak{Com.rings}$ which is monadic and the
cohomology with respect to this monad is known to be
Andr\'{e}-Quillen cohomology. Now we can define a new comonad
$G_{I}(A)(i) = \coprod_{x \rightarrow i} G(A(x))$ on
$\mathfrak{Com.rings}^{I} = \Psi\mathfrak{-rings}$. Using the
bicomplex $C^{*,*}(I,A,M)$ described in the previous section, we can
define cohomology of $\Psi$-rings. If $R: \mathbb{N}^{mult}
\rightarrow \mathfrak{Com.rings}$ is a $\Psi$-ring and $M$ is an
$R$-module, then for any $n \geq 0$, there is a natural system on
$\mathbb{N}^{mult}$ as follows:

$$D_{f} := H_{AQ}^{n}(R,M^{f})$$ where $M^{f}$ is an $R$-module with
$M$ as an abelian group with the following action of $R$:
$$(r,a) \mapsto \Psi^{n}(r)\Psi^{n}(a), \textrm{ for $r \in R, a \in M$}$$
For $u \in \mathbb{N}^{mult}$, we have $u_{*}:D_{f} \rightarrow
D_{uf}$ which is induced by $\Psi^{u}: M^{f} \rightarrow M^{uf}$.
For $v \in \mathbb{N}^{mult}$, we have $v^{*}:D_{f} \rightarrow
D_{fv}$ which is induced by $\Psi^{v}: R \rightarrow R$.

There exists a spectral sequence:
$$E^{p,q}_{2} = H^{p}_{BW}(\mathbb{N}^{mult},H_{AQ}^{q}(A,M)) \Rightarrow H^{p+q}_{\Psi}(A,M).$$
where $H^{*}_{\Psi}(A,M)$ denotes the cohomology of $\Psi$-rings as
it is defined via comonads.

\subsection{$\Pi$-algebras}
A $\Pi$-algebra is a graded group equipped with the action of
primary homotopy operations modeled on the homotopy groups of a
space. Dwyer and Kan \cite{bib:DK} defined the Quillen cohomology of
$\Pi$-algebras, which we denote by $H^{*}_{DK}(A,M)$. Blanc,
Johnson, and Turner \cite{bib:blanc}, defined the Quillen cohomology
of diagrams of $\Pi$-algebras, which we denote by
$H^{*}_{BJT}(A,M)$. However, it is known that Quillen's and Beck's
approaches yield the same cohomology.

An application of our main result is that there exists a spectral
sequence:
$$E^{p,q}_{2} = H^{p}_{BW}(I,\mathcal{H}_{DK}^{q}(A,M)) \Rightarrow H^{p+q}_{BJT}(A,M).$$
where $\mathcal{H}^{*}_{BJT}(A,M)$ is the natural system on $I$
whose value on $\alpha: i \rightarrow j$ is given by
$H^{*}_{DK}(A(i),\alpha^{*}M(j))$.

If we let $I$ be the small category with two distinct objects $0,1$
and one non-trivial map $0 \rightarrow 1$, then our spectral
sequence yields corollary 4.27 in \cite{bib:blanc}.

\subsection{Diagrams of groups}
In the paper  by Cegarra \cite{bib:Cegarra}, the cohomology of
diagrams of groups is described, which we denote by
$H^{*}_{C}(G,A)$. There is also described the following spectral
sequence:

Let $I$ be a small category. If $G: I \rightarrow \mathbf{Gp}$ is an
$I$-group and $A$ is a $G$-module, then for any $n \geq 0$, there is
a natural system on $I$ as follows:

$$\mathcal{H}^{n}(G,A): \mathcal{F}I \rightarrow \mathbf{Ab}, \textrm{\qquad} \xymatrix{u \ar[r]^{\sigma}&v} \mapsto \left\{%
\begin{array}{ll}
    H^{n}(G(u),A(v)) & \hbox{if $n\geq 2$;} \\
    Der(G(u),A(v)) & \hbox{if $n=1$;} \\
    0, & \hbox{if $n=0$.} \\
\end{array}%
\right.    $$ Then there is a natural spectral sequence:
$$E^{p,q}_{2} = H^{p}_{BW}(I, \mathcal{H}^{q+1}(G,A)) \Longrightarrow
H^{p+q+1}_{C}(G,A)$$ where $\mathcal{H}^{q+1}(G,A))$ is the natural
system on $I$ as described above.

\end{document}